\documentclass[a4paper]{amsart}
\usepackage{amsmath,amsthm,amssymb,latexsym,epic,bbm,comment,mathbbol}
\usepackage{graphicx,enumerate,stmaryrd,color}
\usepackage[all,2cell]{xy}
\xyoption{2cell}

\newtheorem{theorem}{Theorem}
\newtheorem{lemma}[theorem]{Lemma}

\newtheorem{proposition}[theorem]{Proposition}

\newtheorem{question}[theorem]{Question}

\newtheorem{problem}[theorem]{Problem}
\usepackage[all]{xy}
\usepackage[active]{srcltx}
\usepackage[parfill]{parskip}
\usepackage{enumerate}

\font\sc=rsfs10
\newcommand{\cC}{\sc\mbox{C}\hspace{1.0pt}}
\newcommand{\cG}{\sc\mbox{G}\hspace{1.0pt}}

\newcommand{\cI}{\sc\mbox{I}\hspace{1.0pt}}
\newcommand{\cJ}{\sc\mbox{J}\hspace{1.0pt}}
\newcommand{\cS}{\sc\mbox{S}\hspace{1.0pt}}

\newcommand{\cA}{\sc\mbox{A}\hspace{1.0pt}}

\newcommand{\cU}{\sc\mbox{U}\hspace{1.0pt}}
\newcommand{\cQ}{\sc\mbox{Q}\hspace{1.0pt}}
\newcommand{\cF}{\sc\mbox{F}\hspace{1.0pt}}

\font\scc=rsfs7
\newcommand{\ccC}{\scc\mbox{C}\hspace{1.0pt}}

\begin{document}

\title[Classification problems in $2$-representation theory]
{Classification problems in\\ $2$-representation theory}
\author{Volodymyr Mazorchuk}

\begin{abstract}
This article surveys  recent advances and future challenges
in the $2$-representation theory of finitary $2$-categories with 
a particular emphasis on problems related to classification of 
various classes of $2$-representations.
\end{abstract}

\maketitle

\section{Introduction}\label{s1}

Higher representation theory emerged as an offspring of categorification.
The latter term traditionally describes the approach, originated in 
\cite{CF,Cr}, of upgrading set-theoretical notions to category theoretical with
a hope to create more structure. The major breakthrough of categorification
was invention of Khovanov homology in \cite{Kh}. After that several other 
spectacular applications followed, for example, to Brou{\'e}'s abelian defect
group conjecture in \cite{CR}, to study of modules over Lie superalgebras
in \cite{BS,BLW} and to various other problems. The common feature of all these 
and many other applications is that construction and comparison of functorial 
actions on different categories was the key part of the argument.

Functorial actions form the red thread of \cite{BFK} leading to an alternative
reformulation of Khovanov homology in \cite{St} using BGG category $\mathcal{O}$. 
In \cite{CR} and, later on, in \cite{Ro}, functorial actions were abstractly 
reformulated in terms of representation theory of certain $2$-categories. 
This direction of study was subsequently  called {\em higher representation
theory} or, alternatively, just {\em $2$-representation theory} to emphasize that, so far,
it directs only at this second level of general higher categories. However, one has to
note that these and many further papers like \cite{KL,El} and other mainly study 
special examples of $2$-categories which originate in topologically motivated 
diagrammatic calculus.

The series \cite{MM1,MM2,MM3,MM4,MM5,MM6} of papers took higher representation theory
to a more abstract level. These papers started a systematic study of abstract
$2$-analogues of finite dimensional algebras, called {\em finitary $2$-categories},
and their representation theory. Despite of the fact that the $2$-category of 
$2$-representations of an abstract finitary $2$-category is much more complicated
than the category of modules over a finite dimensional algebra (in particular,
it is not abelian), it turned out
that, in many cases, it is possible to {\em construct}, {\em compare} and even 
{\em classify} various classes of $2$-representations. These kinds of problems
were recently studied in a number of papers, see \cite{MM5,MM6,Zh2,MZ1,MaMa,KMMZ,MT,MMMT,MMZ,MZ2}.
The aim of the present article is to give an overview of these results with a particular 
emphasis on open problems and future challenges.

We start in Section~\ref{s2} with a brief description of main objects of study.
Section~\ref{s3} lists a number of classical examples. Sections~\ref{s4},
\ref{s5} and \ref{s6} then concentrate on special classes of $2$-representations.
In particular, Sections~\ref{s4} addresses {\em cell $2$-representations}
and Sections~\ref{s5} is devoted to {\em simple transitive $2$-representations}.
We do not provide any proofs but rather give explicit references to original
sources.

The present paper might serve as a complement reading to the series of lectures on
higher representation theory which the author gave during 
{\em Brazilian Algebra Meeting} in Diamantina, Brazil, in August 2016.
\vspace{5mm}

{\bf Acknowledgments.} The author is partially supported by the
Swedish Research Council and G{\"o}ran Gustafsson Stiftelse.
The authors thanks the organizers of the Brazilian Algebra Meeting
for invitation to give the series of lecture on higher representation theory.

\section{Finitary $2$-categories and their $2$-representations}\label{s2}

\subsection{Finitary $2$-categories}\label{s2.1}

A {\em $2$-category} is a category which is enriched over the monoidal category
$\mathbf{Cat}$ of small categories (the monoidal structure of the latter category is given by
Cartesian product). This means that a $2$-category $\cC$ consists of 
\begin{itemize}
\item objects, denoted $\mathtt{i}$, $\mathtt{j}$, $\mathtt{k}$,\dots;
\item small morphism categories $\cC(\mathtt{i},\mathtt{j})$, for all $\mathtt{i},\mathtt{j}\in\cC$;
\item bifunctorial compositions $\cC(\mathtt{j},\mathtt{k})\times \cC(\mathtt{i},\mathtt{j})
\to \cC(\mathtt{i},\mathtt{k})$, for all $\mathtt{i},\mathtt{j},\mathtt{k}\in\cC$;
\item and identity objects $\mathbbm{1}_{\mathtt{i}}\in \cC(\mathtt{i},\mathtt{i})$, 
for all $\mathtt{i}\in\cC$;
\end{itemize}
which satisfy the obvious collection of strict axioms. The following terminology is standard:
\begin{itemize}
\item objects in $\cC(\mathtt{i},\mathtt{j})$ are called {\em $1$-morphisms} of $\cC$
and will be denoted by $\mathrm{F}$, $\mathrm{G}$, etc.;
\item morphisms in $\cC(\mathtt{i},\mathtt{j})$ are called {\em $2$-morphisms} of $\cC$
and will be denoted by $\alpha$, $\beta$, etc.;
\item composition of $2$-morphisms inside $\cC(\mathtt{i},\mathtt{j})$ is called {\em vertical}
and will be denoted by $\circ_v$;
\item composition of $2$-morphisms coming from
$\cC(\mathtt{j},\mathtt{k})\times \cC(\mathtt{i},\mathtt{j})
\to \cC(\mathtt{i},\mathtt{k})$ is called {\em horizontal}
and will be denoted by $\circ_h$. 
\end{itemize}
As usual, for a $1$-morphism $\mathrm{F}$, the identity $2$-morphism for $\mathrm{F}$ 
is denoted $\mathrm{id}_{\mathrm{F}}$, moreover, for a $2$-morphism $\alpha$, 
the compositions $\mathrm{id}_{\mathrm{F}}\circ_h\alpha$ and $\mathrm{id}_{\mathrm{F}}\circ_h\alpha$
are denoted by $\mathrm{F}(\alpha)$ and $\alpha_{\mathrm{F}}$, respectively.
We refer the reader to \cite{Le,Mc} for more details on $2$-categories and for various generalizations,
in particular, for the corresponding {\em non-strict} notion of a {\em bicategory}.

Let $\Bbbk$ be an algebraically closed field. Recall that a category $\mathcal{C}$ is called
{\em finitary $\Bbbk$-linear} provided that it is equivalent to the category of projective modules
over a finite dimensional (associative) $\Bbbk$-algebra. Each such category is $\Bbbk$-linear,
that is enriched over the category of $\Bbbk$-vector spaces, moreover, it is idempotent split and 
Krull-Schmidt and has finitely many isomorphism classes of indecomposable
objects and finite dimensional (over $\Bbbk$) spaces of morphisms.

Following \cite[Subsection~2.2]{MM1}, we will say that a {\em $2$-category} $\cC$ is 
{\em finitary over $\Bbbk$} provided that
\begin{itemize}
\item each  $\cC(\mathtt{i},\mathtt{j})$ is finitary $\Bbbk$-linear;
\item all compositions are biadditive and $\Bbbk$-bilinear, whenever appropriate;
\item all $\mathbbm{1}_{\mathtt{i}}$ are indecomposable.
\end{itemize}
The last condition is a technical condition which makes the life easier at many occasions. From the representation
theoretic prospective, this condition is not restrictive as, starting from a $2$-category satisfying all
other conditions, one can use idempotent splitting to produce a finitary $2$-category with essentially 
the same representation theory.

In what follows, we will simply say that $\cC$ is {\em finitary} as our field $\Bbbk$ will be
fixed throughout the paper (with the exception of examples related to Soergel bimodules where
$\Bbbk=\mathbb{C}$).

\subsection{$2$-representations}\label{s2.2}

For two $2$-categories $\cA$ and $\cC$, a {\em $2$-functor} $\Phi:\cA\to \cC$ is a functor which 
respects all $2$-categorical structure. This means that $\Phi$ 
\begin{itemize}
\item maps $1$-morphisms to $1$-morphisms;
\item maps $2$-morphisms to $2$-morphisms;
\item is compatible with composition of $1$-morphisms;
\item is compatible with both horizontal and vertical composition of $2$-morphisms;
\item sends identity $1$-morphisms to identity $1$-morphisms;
\item sends identity $2$-morphisms to identity $2$-morphisms.
\end{itemize}

A {\em $2$-representation} of a $2$-category $\cC$ is a $2$-functor to some fixed $2$-category.
Classical examples of such target $2$-categories are:
\begin{itemize}
\item the $2$-category $\mathbf{Cat}$ of small categories, here  $1$-morphisms are functors and $2$-morphisms
are natural transformations;
\item the $2$-category $\mathfrak{A}^{f}_{\Bbbk}$ of {\em finitary $\Bbbk$-linear} categories, here
$1$-morphisms are additive $\Bbbk$-linear functors and $2$-morphisms are natural transformations;
\item the $2$-category $\mathfrak{R}_{\Bbbk}$ of {\em finitary abelian $\Bbbk$-linear} categories, here
objects are categories equivalent to module categories of finite dimensional (associative) $\Bbbk$-algebras,
$1$-morphisms are right exact additive $\Bbbk$-linear functors and $2$-morphisms are natural transformations.
\end{itemize}

All $2$-representations of a $2$-category $\cC$ (in a fixed target $2$-category) form a $2$-category. 
In this $2$-category we have:
\begin{itemize}
\item $1$-morphisms are (strong) $2$-natural transformations;
\item $2$-morphisms are modifications.
\end{itemize}
One has to make a choice for the level of strictness for $1$-morphisms in the 
$2$-category of $2$-representations of $\cC$. In the language of \cite[Subsection~1.2]{Le}, this corresponds
to choosing between the so-called strong or strict transformations. Strict transformations were considered
in the paper \cite{MM1} and the setup was changed to strong transformations in \cite{MM3}. The latter
allows for more flexibility and more reasonable results (for example, the relation of equivalence of 
two $2$-representations becomes symmetric).

For a finitary $2$-category $\cC$, its $2$-representations in $\mathfrak{A}^{f}_{\Bbbk}$ are called
{\em finitary additive} $2$-representations and the corresponding $2$-category of $2$-representations
is denoted by $\cC\text{-}\mathrm{afmod}$. Further, $2$-representations of $\cC$ in $\mathfrak{R}_{\Bbbk}$ 
are called {\em abelian} $2$-representations and the corresponding $2$-category of $2$-representations
is denoted by $\cC\text{-}\mathrm{mod}$. Note that neither $\cC\text{-}\mathrm{afmod}$ nor
$\cC\text{-}\mathrm{mod}$ are abelian categories.

We will usually denote $2$-representations of $\cC$ by $\mathbf{M}$, $\mathbf{N}$ etc. For the
sake of readability, it is often convenient to use the actions notation $\mathrm{F}\, X$ instead
of the representation notation $\mathbf{M}(\mathrm{F})\big( X\big)$.

Here is an example of a $2$-representation: for $\mathtt{i}\in\cC$, the {\em principal} 
$2$-representation $\mathbf{P}_{\mathtt{i}}$ is defined to be the Yoneda $2$-representation
$\cC(\mathtt{i},{}_-)$. If $\cC$ is finitary, we have $\mathbf{P}_{\mathtt{i}}\in \cC\text{-}\mathrm{afmod}$.

Two $2$-representations $\mathbf{M}$ and $\mathbf{N}$ of $\cC$ are called {\em equivalent} provided 
that there is a $2$-natural transformation $\Phi:\mathbf{M}\to\mathbf{N}$ such that 
$\Phi_{\mathtt{i}}$ is an equivalence of categories, for each $\mathtt{i}\in\cC$.

\subsection{Abelianization}\label{s2.3}

Given a finitary $\Bbbk$-linear category $\mathcal{C}$, the {\em diagrammatic abe\-li\-a\-ni\-za\-tion} of
$\mathcal{C}$ is the category $\overline{\mathcal{C}}$ of diagrams of the form 
$X\overset{\alpha}{\longrightarrow} Y$ over $\mathcal{C}$ with morphisms being the obvious commutative
squares modulo the projective homotopy relations. The category $\overline{\mathcal{C}}$ is abelian
and is equivalent to the category of modules over a finite dimensional $\Bbbk$-algebra.
The original category $\mathcal{C}$ canonically embeds into $\overline{\mathcal{C}}$ via diagrams
of the form $0\to Y$ and this embedding provides an equivalence between $\mathcal{C}$ and the
category of projective objects in $\overline{\mathcal{C}}$. We refer the reader to \cite{Fr} for details.

For a finitary $2$-category $\cC$, using diagrammatic abelianization and component-wise action on 
diagrams defines a $2$-functor 
\begin{displaymath}
\overline{\hspace{1mm}\cdot\hspace{1mm}}: \cC\text{-}\mathrm{afmod}\to \cC\text{-}\mathrm{mod},
\end{displaymath}
called {\em abelianization}, see \cite[Subsection~3.1]{MM1}. In \cite[Section~3]{MMMT} one finds
a more advanced refinement of this construction which is way more technical but also 
has some extra nice properties.

\subsection{Fiat $2$-categories}\label{s2.4}

As we will see, many examples of finitary $2$-categories have additional structure which plays
very important role and significantly simplifies many arguments. This additional structure, on a finitary
$2$-category $\cC$, consists of
\begin{itemize}
\item a {\em weak involution} $\star$ which inverts the direction of both $1$- and $2$-morphisms,
\item {\em adjunction morphisms} 
$\varepsilon^{(\mathrm{F})}:\mathrm{F}\circ \mathrm{F}^{\star}\to\mathbbm{1}_{\mathtt{j}}$
and $\eta^{(\mathrm{F})}:\mathbbm{1}_{\mathtt{i}}\to \mathrm{F}^{\star}\circ \mathrm{F}$,
for each $\mathrm{F}\in\cC(\mathtt{i},\mathtt{j})$, which make $(\mathrm{F},\mathrm{F}^{\star})$
a pair of adjoint $1$-morphisms in the sense that 
\begin{displaymath}
\mathrm{id}_{\mathrm{F}}=\varepsilon^{(\mathrm{F})}_{\mathrm{F}}\circ_v\mathrm{F}(\eta^{(\mathrm{F})})
\quad\text{ and }\quad
\mathrm{id}_{\mathrm{F}}^{\star}=\mathrm{F}^{\star}(\varepsilon^{(\mathrm{F})})\circ_v
\eta^{(\mathrm{F})}_{\mathrm{F}^{\star}}.
\end{displaymath}
\end{itemize}
A $2$-category $\cC$ having such an additional structure is called {\em fiat}, where
``f'' stands for {\em finitary}, ``i'' stands for {\em involution}, ``a'' stands for 
{\em adjunction} and  ``t'' stands for {\em $2$-category}, see \cite[Subsection~2.4]{MM1}. 
If a similar structure exists for a not necessarily involutive anti-autoequivalence $\star$,
the $2$-category $\cC$ is called {\em weakly fiat}, see \cite[Subsection~7.3]{MM2} and \cite[Appendix]{MM6}.

In many situations, involutions in $2$-categories change the direction of $1$-morphisms but preserve the
direction of $2$-morphisms, see e.g. \cite[Page~3]{Le}. The above definition, in which both the
directions of $1$- and $2$-morphisms get reversed, is motivated by the
$2$-category of endofunctors of $A$-mod, for a finite dimensional $\Bbbk$-algebra $A$. For each
pair $(\mathrm{F},\mathrm{G})$ of adjoint endofunctors of $A$-mod, there is an $A$-$A$-bimodule $Q$ such that
$\mathrm{F}$ is isomorphic to $Q\otimes_A{}_-$ and $\mathrm{G}$ is isomorphic to 
$\mathrm{Hom}_{A-}(Q,{}_-)$, see \cite[Chapter~I]{Ba}. Natural transformations between functors 
correspond to homomorphisms between bimodules. When taking the adjoint functor, the bimodule $Q$ 
ends up on the contravariant place of the bifunctor $\mathrm{Hom}$ and hence the direction of natural
transformations gets reversed.

In the literature one could find similar structures under the name of {\em rigid} 
tensor categories categories, see e.g. \cite{EGNO}. 

\subsection{Grothendieck decategorification}\label{s2.5}

For an finitary $\Bbbk$-linear category $\mathcal{C}$, let $[\mathcal{C}]_{\oplus}$ denote the
{\em split Grothendieck group} of $\mathcal{C}$. Then $[\mathcal{C}]_{\oplus}$ is a free abelian group
which has a canonical  generating set given by isomorphism classes of indecomposable objects in $\mathcal{C}$.
For a category $\mathcal{C}$ equivalent to $A$-mod, for some finite dimensional $\Bbbk$-algebra $A$,
let $[\mathcal{C}]$ denote the {\em Grothendieck group} of $\mathcal{C}$. Then $[\mathcal{C}]$ is a free 
abelian group which has a canonical  generating set given by isomorphism classes of simple objects in $\mathcal{C}$.

Let $\cC$ be a finitary $2$-category and $\mathbf{M}\in \cC$-afmod. Let $[\cC]_{\oplus}$ denote the
ordinary category which has  the same objects as $\cC$ and in which morphisms are given by 
$[\cC]_{\oplus}(\mathtt{i},\mathtt{j}):=[\cC(\mathtt{i},\mathtt{j})]_{\oplus}$ with induced composition.
The category $[\cC]_{\oplus}$ is called the {\em Grothendieck decategorification} of $\cC$. The category
$[\cC]_{\oplus}$ acts on abelian groups $[\mathbf{M}(\mathtt{i})]_{\oplus}$, where $\mathtt{i}\in \cC$,
which in this way defines the {\em Grothendieck decategorification} $[\mathbf{M}]_{\oplus}$ of $\mathbf{M}$.

If $\cC$ is fiat and $\mathbf{M}\in \cC$-mod, then $[\cC]_{\oplus}$ acts on abelian groups 
$[\mathbf{M}(\mathtt{i})]$, where $\mathtt{i}\in \cC$,
which in this way defines the {\em Grothendieck decategorification} $[\mathbf{M}]$ of $\mathbf{M}$.

We note that there are alternative decategorifications, notably, the {\em trace decategorification}
introduced in \cite{BGHL}.

\section{Examples of finitary $2$-categories}\label{s3}

\subsection{Set-theoretic issues}\label{s3.1}

There are some set-theoretic complications due to the fact that, by definition,
each $\cC(\mathtt{i},\mathtt{j})$ of a $2$-category $\cC$ has to be small. This
prevents us to consider, for example, the category of all $A$-$A$--bimodules, for
a finite dimensional $\Bbbk$-algebra $A$, as a $2$-category (with one formal object
that can be identified with $A$-mod). The reason for that is the observation
that the category of all bimodules is not small. 

In what follows we will give many examples by considering all endofunctors (of some type)
of some category $\mathcal{C}$. To avoid the above problem, we will always assume that 
$\mathcal{C}$ is small. This, however, creates a choice. For example, one has to choose
a small category $\mathcal{C}$ equivalent to $A$-mod. Different choices of $\mathcal{C}$
lead to different, however, usually {\em biequivalent}, $2$-categories.

\subsection{Projective endofunctors of $A$-mod}\label{s3.2}

Let $A$ be a finite dimensional, basic,  connected $\Bbbk$-algebra. Fix a small category 
$\mathcal{C}$ equivalent to $A$-mod. Recall that a {\em projective} $A$-$A$-bimodule
is an $A$-$A$-bimodule from the additive closure $\mathrm{add}(A\otimes_{\Bbbk}A)$
of ${}_AA\otimes_{\Bbbk}A_A$. We also have the {\em regular} or {\em identity} $A$-$A$-bimodule 
${}_AA_A$.

Denote by $\cC_A=\cC_{A,\mathcal{C}}$ the $2$-category which has
\begin{itemize}
\item one object $\mathtt{i}$ (which should be though of as $\mathcal{C}$);
\item as $1$-morphisms, endofunctors of $\mathcal{C}$ from the additive closure of 
endofunctors given by tensoring with projective or regular $A$-$A$-bimodules;
\item as $2$-morphisms, all natural transformations of functors.
\end{itemize}
This $2$-category appears in \cite[Subsection~7.3]{MM1}.

If $A$ is simple, then $A\cong\Bbbk$ and $\cC_A$ has a unique (up to isomorphism)
indecomposable $1$-morphism, namely $\mathbb{1}_{\mathtt{i}}$. 
If $A$ is not simple, let $e_1+e_2+\dots+e_n=1$ be a primitive decomposition of the
identity $1\in A$. Then, apart from $\mathbb{1}_{\mathtt{i}}$, the $2$-category $\cC_A$ has
$n^2$ additional indecomposable $1$-morphisms $\mathrm{F}_{ij}$, where,
for $i,j=1,2,\dots,n$, the morphism $\mathrm{F}_{ij}$ corresponds to tensoring
with $Ae_i\otimes_{\Bbbk}e_jA$. We note that 
\begin{displaymath}
\mathrm{F}_{ij}\circ\mathrm{F}_{st}\cong \mathrm{F}_{it}^{\dim(e_jAe_s)}\quad\text{ and }\quad
\mathrm{F}\circ \mathrm{F}\cong\mathrm{F}^{\dim(A)},\quad\text{ for }\quad
\mathrm{F}:=\bigoplus_{i,j=1}^n\mathrm{F}_{ij}.
\end{displaymath}
The $2$-category $\cC_A$ is always finitary. It is weakly fiat if and only if $A$ is self-injective.
It is fiat if and only if $A$ is weakly symmetric. In the latter case, 
$(\mathrm{F}_{ij},\mathrm{F}_{ji})$ forms an adjoint pair of $1$-morphisms, for all $i$ and $j$.

\subsection{Finitary $2$-categories of all bimodules}\label{s3.3}

For $n=1,2,\dots$, let $A_n$ denote the $\Bbbk$-algebra given as the quotient of the
path algebra of the quiver
\begin{equation}\label{eq1}
\xymatrix{1\ar[rr]&&2\ar[rr]&&3\ar[rr]&&\dots\ar[rr]&&n} 
\end{equation}
by the relations that the product of any two arrows is zero. Let $\mathcal{C}$ be a small
category equivalent to $A$-mod. 

Denote by $\cF_{A_n}=\cF_{A_n,\mathcal{C}}$ the $2$-category which has
\begin{itemize}
\item one object $\mathtt{i}$ (which should be though of as $\mathcal{C}$);
\item as $1$-morphisms, all right exact endofunctors of $\mathcal{C}$;
\item as $2$-morphisms, all natural transformations of functors.
\end{itemize}
The $2$-category $\cF_{A_n}$ is finitary, see \cite[Section~2]{MZ2}. The reason for that is the fact that
the enveloping algebra $A_n\otimes_{\Bbbk}A_n^{\mathrm{op}}$ of $A_n$ is a special biserial algebra
in the sense of \cite{BR,WW}, so one can use the classification of indecomposable modules for this
algebra to check that it has, in fact, only finitely many of them, up to isomorphism. Unless $n=1$,
the $2$-category $\cF_{A_n}$ is neither fiat nor weakly fiat.

Note that a similar $2$-category $\cF_{A}$ can be defined for any connected and basic $\Bbbk$-algebra $A$.
However, if $\cF_{A}$  is fiat, then $A\cong A_n$, for some $n$, see \cite[Theorem~1]{MZ2}.

\subsection{Subbimodules of the identity bimodule}\label{s3.4}

Let $\Gamma$ be a finite oriented tree and $A$ be the path algebra of $\Gamma$.
Let $\mathcal{C}$ be a small category equivalent to $A$-mod. 

Denote by $\cG_{A}=\cG_{A,\mathcal{C}}$ the $2$-category which has
\begin{itemize}
\item one object $\mathtt{i}$ (which should be though of as $\mathcal{C}$);
\item as $1$-morphisms, all endofunctors of $\mathcal{C}$ from the additive closure of
subfunctors of the identity functor;
\item as $2$-morphisms, all natural transformations of functors.
\end{itemize}
The $2$-category $\cG_{A}$ is finitary due to the fact that the regular $A$-$A$-bimodule is
multiplicity free and hence has only finitely many subbimodules. Unless $\Gamma$ has one vertex,
the $2$-category $\cF_{A_n}$ is neither fiat nor weakly fiat.

The $2$-category $\cG_{A}$ first appeared, slightly disguised, in \cite{GM1} 
(inspired by \cite{Gr}), for $\Gamma$ being the quiver in \eqref{eq1}. 
It was further studied, for various special types of trees, in \cite{GM2} and \cite{Zh1,Zh2}.

\subsection{Soergel bimodules for finite Coxeter systems}\label{s3.5}

Let $(W,S)$ be a finite Coxeter system and $\mathfrak{h}$ be a reflection faithful complexified $W$-module.
Let, further, $\mathtt{C}$ be the corresponding coinvariant algebra, that is the quotient of 
$\mathbb{C}[\mathfrak{h}]$ by the ideal generated by homogeneous $W$-invariant polynomials of positive degree.
Then $\mathtt{C}$ is a finite dimensional algebra which carries the natural structure of a
regular $W$-module, in particular, $\dim(\mathtt{C})=|W|$. 

For $s\in S$, denote by $\mathtt{C}^s$ the subalgebra of $s$-invariants in $\mathtt{C}$. Then
$\mathtt{C}$ is a free $\mathtt{C}^s$-module of rank two. Both algebras $\mathtt{C}$ and 
$\mathtt{C}^s$ are symmetric (recall that $A$ is symmetric provided that ${}_AA_A\cong
\mathrm{Hom}_{\Bbbk}({}_AA_A,\Bbbk)$). We refer to \cite{Hi} for details.

For $w\in W$ with reduced decomposition $w=s_1s_2\cdots s_k$, define the {\em Bott-Samelson}
$\mathtt{C}$-$\mathtt{C}$-bimodule 
\begin{displaymath}
\hat{B}_w:=\mathtt{C}\otimes_{\mathtt{C}^{s_1}}\mathtt{C}\otimes_{\mathtt{C}^{s_2}}\mathtt{C}
\otimes_{\mathtt{C}^{s_3}}\dots \otimes_{\mathtt{C}^{s_k}}\mathtt{C}. 
\end{displaymath}

Let $\mathcal{C}$ be a small category equivalent to $\mathtt{C}$-mod. 
Denote by $\cS_{W}=\cS_{W,S,V,\mathcal{C}}$ the $2$-category which has
\begin{itemize}
\item one object $\mathtt{i}$ (which should be though of as $\mathcal{C}$);
\item as $1$-morphisms, all endofunctors of $\mathcal{C}$ coming from tensoring with bimodules
from the additive closure of all Bott-Samelson bimodules;
\item as $2$-morphisms, all natural transformations of functors.
\end{itemize}
The $2$-category $\cS_{W}$ is called the $2$-category of {\em Soergel bimodules} over $\mathtt{C}$.
The non-trivial point here is the fact that it is closed under composition of functors. This is
shown in \cite{So2}. Moreover, it is also shown in \cite{So2} that, for each $w\in W$, the 
bimodule $\hat{B}_w$ contains a unique indecomposable summand $B_w$ which does not belong to the 
additive closure of all $\hat{B}_x$, where the length of $x$ is strictly smaller than that of $w$. 
Bimodules $B_w$ are usually called {\em Soergel bimodules} (although the name is also used
for any direct sum of such bimodules). The $2$-category $\cS_{W}$ is both
finitary and fiat, where $(B_w,B_{w^{-1}})$ forms an adjoint pair of $1$-morphisms, for all $w$.

The theory is inspired by \cite{So1} where, for finite {\em Weyl} groups, Soergel bimodules
appear as ``combinatorial description'' of indecomposable projective functors on the principal
block of the BGG category $\mathcal{O}$ associated with a triangular decomposition of the 
simple finite dimensional Lie algebra corresponding to $W$, see \cite{BGG,BG,Hu} for details on the latter.
An explicit connection between the $2$-category of Soergel bimodules and the
Kazhdan-Lusztig basis of the  Hecke algebra of $(W,S)$ was established in \cite{EW}.

\subsection{Singular Soergel bimodules}\label{s3.6}

Let $(W,S)$ be a finite Coxeter system, $\mathfrak{h}$ a reflection faithful complexified $W$-module and
$\mathtt{C}$ the corresponding coinvariant algebra. For each $T\subset S$, let 
$W^T$ be the subgroup of $W$ generated by all $t\in T$ and $\mathtt{C}^T$ the subalgebra of
$W^T$-invariants in $\mathtt{C}$. Note that $\mathtt{C}^{\varnothing}=\mathtt{C}$ while
$\mathtt{C}^S=\mathbb{C}$.

For each $T\subset S$, let $\mathcal{C}^T$ be a small category equivalent to $\mathtt{C}^T$-mod. 
Denote by $\cS\cS_{W}=\cS\cS_{W,S,V,\mathcal{C}}$ the $2$-category which has
\begin{itemize}
\item objects $\mathtt{i}^T$, where $T\subset S$ (each $\mathtt{i}^T$ should be 
though of as $\mathcal{C}^T$);
\item as $1$-morphisms, all endofunctors from  $\mathcal{C}^T$ to $\mathcal{C}^R$ 
given by 
\begin{displaymath}
\mathrm{Res}^{\mathtt{C}}_{\mathtt{C}^R} \circ\mathrm{F}\circ \mathrm{Ind}^{\mathtt{C}}_{\mathtt{C}^T},
\end{displaymath}
where $\mathrm{F}$ is given by a usual Soergel bimodule;
\item as $2$-morphisms, all natural transformations of functors.
\end{itemize}
The $2$-category $\cS\cS_{W}$ is called the $2$-category of {\em singular Soergel bimodules} over $\mathtt{C}$.
The $2$-category $\cS\cS_{W}$ is both finitary and fiat.

The $2$-category $\cS\cS_{W}$ also admits an alternative description using projective
functors between {\em singular} blocks of the BGG category $\mathcal{O}$, see \cite{So1,BGG,BG}.

\subsection{Finitary quotients of $2$-Kac-Moody algebras}\label{s3.7}

Let $\mathfrak{g}$ be a simple complex finite dimensional Lie algebra and $\dot{U}_{\mathfrak{g}}$ the
idempotented version of the universal enveloping algebra of $\mathfrak{g}$, see \cite{Lu}.
The papers \cite{KL,Ro} introduce certain (not finitary)  $2$-categories whose Grothendieck decategorification
is isomorphic to the integral form of $\dot{U}_{\mathfrak{g}}$. In \cite{Br} it is further shown that the
two (slightly different) constructions in \cite{KL,Ro} give, in fact, biequivalent $2$-categories.
These are the so-called {\em $2$-Kac-Moody algebras} of finite type which we will denote $\cU_{\mathfrak{g}}$.

Each simple finite dimensional $\mathfrak{g}$-module $V(\lambda)$, where $\lambda$ is the highest weight, 
admits a {\em categorification} in the sense that there exists a $2$-representation $\mathbf{M}_{\lambda}$ 
of $\cU_{\mathfrak{g}}$ (even a unique one, up to equivalence, under the additional assumption that the 
object $\lambda$ is represented by a non-zero semi-simple category) whose Grothendieck decategorification 
is isomorphic to the integral form of $V(\lambda)$. Let $\cI_{\lambda}$ be the kernel of $\mathbf{M}_{\lambda}$.
Then $\cI_{\lambda}$ is a two-sided $2$-ideal of $\cU_{\mathfrak{g}}$ and the quotient $2$-category
$\cU_{\mathfrak{g}}/\cI_{\lambda}$ is both, finitary and fiat.

\section{Cells and cell $2$-representations}\label{s4}

\subsection{Cells}\label{s4.1}

For a finitary $2$-category  $\cC$, denote by $\mathcal{S}[\cC]$ the set of isomorphism classes of indecomposable
$1$-objects in $\cC$. The set $\mathcal{S}[\cC]$ is finite and has the natural multivalued operation
$\bullet$ given, for $[\mathrm{F}],[\mathrm{G}]\in \mathcal{S}[\cC]$, by 
\begin{displaymath}
[\mathrm{F}]\bullet[\mathrm{G}]=
\{[\mathrm{H}]\,:\,\mathrm{H}\text{ is isomorphic to a summand of } \mathrm{F}\circ \mathrm{G}\}.
\end{displaymath}
The operation $\bullet$ is associative (as a  multivalued operation) and hence defines on $\mathcal{S}[\cC]$
the structure of a {\em multisemigroup}, see \cite[Section~3]{MM2} (see also \cite{KuMa} for more details on multisemigroups).

The {\em left} partial pre-order $\leq_L$ on $\mathcal{S}[\cC]$ is defined by setting 
$[\mathrm{F}]\leq_L[\mathrm{G}]$, for $[\mathrm{F}],[\mathrm{G}]\in \mathcal{S}[\cC]$, provided 
that $\mathrm{G}$ is isomorphic to a summand of $\mathrm{H}\circ \mathrm{F}$, for some $1$-morphism
$\mathrm{H}$. Equivalence classes with respect to $\leq_L$ are called {\em left cells} and the corresponding
equivalence relation is denoted $\sim_L$. 

The {\em right} partial pre-order $\leq_R$, the {\em right cells} and the corresponding
equivalence relation $\sim_R$ are defined similarly using multiplication with  $\mathrm{H}$ on the right 
of $\mathrm{F}$. The {\em two-sided} partial pre-order $\leq_J$, the {\em two-sided cells} and the corresponding
equivalence relation $\sim_J$ are defined similarly using multiplication with  $\mathrm{H}_1$ and
$\mathrm{H}_2$ on both sides of $\mathrm{F}$. 

These notions are similar and spirit to and generalize the notions of {\em Green's} relations and partial
orders for semigroups, see \cite{Gre}, and also the notions of Kazhdan-Lusztig cells and order in
\cite{KaLu}, see also \cite{KM2}. We refer the reader to \cite{KuMa} for more details and to
\cite{KM2} for a generalization to positively based algebras.

For simplicity, we will say ``left cells in $\cC$'' instead of ``left cells in $\mathcal{S}[\cC]$''
and similarly for right and $2$-sided cells.

A two-sided cell $\mathcal{J}$ is said to be {\em regular} provided that 
\begin{itemize}
\item any pair of left cells inside $\mathcal{J}$ is incomparable with respect to the left order;
\item any pair of right cells inside $\mathcal{J}$ is incomparable with respect to the right order.
\end{itemize}
A two-sided cell $\mathcal{J}$ is said to be {\em strongly regular} provided that 
it is regular and $|\mathcal{L}\cap\mathcal{R}|=1$, for any left cell $\mathcal{L}$ in $\mathcal{J}$
and any right cell $\mathcal{R}$ in $\mathcal{J}$. We refer to \cite[Subsection~4.8]{MM1} for details.
A two-sided cell $\mathcal{J}$ is said to be {\em idempotent} provided that it contains three
elements $\mathrm{F}$, $\mathrm{G}$ and  $\mathrm{H}$ (not necessarily distinct) such that 
$\mathrm{F}$ is isomorphic to a direct summand of $\mathrm{G}\circ \mathrm{H}$,
see \cite[Subsection~2.3]{CM}. The following is \cite[Corollary~19]{KM2}.

\begin{proposition}\label{prop5}
Each idempotent two-sided cell of finitary $2$-category is regular, in particular,
each two-sided cell of (weakly) fiat $2$-category is regular.
\end{proposition}

\subsection{Cell $2$-representations}\label{s4.2}

Let $\cC$ be a finitary $2$-category and $\mathcal{L}$ a left cell in $\cC$. Then there is
an object $\mathtt{i}=\mathtt{i}_{\mathcal{L}}\in\cC$ such that all $1$-morphisms in $\mathcal{L}$ 
start from $\mathtt{i}$. Consider the principal $2$-representation $\mathbf{P}_{\mathtt{i}}$.

For each $\mathtt{j}\in\cC$, denote by $\mathbf{M}(\mathtt{j})$ the additive closure in 
$\mathbf{P}_{\mathtt{i}}(\mathtt{j})$ of all indecomposable $1$-morphisms
$\mathrm{F}\in \mathbf{P}_{\mathtt{i}}(\mathtt{j})=\cC(\mathtt{i},\mathtt{j})$ satisfying $\mathcal{L}\leq_L \mathrm{F}$
(note that the latter notation makes sense as $\mathcal{L}$ is a left cell). Then $\mathbf{M}$
has the natural structure of a $2$-representation of $\cC$ which is inherited from 
$\mathbf{P}_{\mathtt{i}}$ by restriction. The following lemma can be found e.g. in \cite[Lemma~3]{MM5}.

\begin{lemma}\label{lem1}
The  $2$-representation $\mathbf{M}$ has a unique maximal 
$\cC$-invariant ideal $\mathbf{I}$. 
\end{lemma}

The quotient $2$-representation $\mathbf{M}/\mathbf{I}$ is called the {\em cell $2$-representation}
corresponding to $\mathcal{L}$ and denoted $\mathbf{C}_{\mathcal{L}}$. The construction presented here
first appears in \cite[Subsection~6.5]{MM2}.
It follows directly from the  construction that isomorphism classes of indecomposable objects in 
\begin{displaymath}
\coprod_{\mathtt{j}\in\ccC} \mathbf{C}_{\mathcal{L}}(\mathtt{j})
\end{displaymath}
correspond bijectively to elements in $\mathcal{L}$. Note that two cell $2$-representations
$\mathbf{C}_{\mathcal{L}}$ and $\mathbf{C}_{\mathcal{L}'}$ might be equivalent even 
in the case $\mathcal{L}\neq \mathcal{L}'$. The following is \cite[Theorem~3.1]{MM6}.

\begin{theorem}\label{thm7}
Let $\cC$ be a weakly fiat $2$-category in which all two-sided cells are strongly regular.
Then, for any left cells $\mathcal{L}$ and $\mathcal{L}'$ in $\cC$, we have 
$\mathbf{C}_{\mathcal{L}}\cong \mathbf{C}_{\mathcal{L}'}$ if and only if 
$\mathcal{L}$ and $\mathcal{L}'$ belong to the same two-sided cell in $\cC$.
\end{theorem}

Cell $2$-representations can be viewed as natural $2$-analogues of semigroups representations
associated to left cells, see \cite[Subsection~11.2]{GM}.

\subsection{Basic properties of cell $2$-representations}\label{s4.3}

Let $\cC$ be a finitary $2$-category and $\mathbf{M}\in \cC$-afmod. We will say that $\mathbf{M}$
is {\em transitive} provided that, for any $\mathtt{i},\mathtt{j}\in \cC$ and for any 
indecomposable $X\in \mathbf{M}(\mathtt{i})$ and $Y\in \mathbf{M}(\mathtt{j})$, there is a
$1$-morphism $\mathrm{F}\in\cC(\mathtt{i},\mathtt{j})$ such that $Y$ is isomorphic to a
summand of $\mathrm{F}\, X$.

Directly from the definition of left cells and the construction of cell $2$-representations
it follows that each cell $2$-representation is transitive.

We will say that $\mathbf{M}$ is {\em simple} provided that $\mathbf{M}$ has no proper 
$\cC$-stable ideals. Note that every simple  $2$-representation is automatically transitive.

Directly from construction of cell $2$-representations and Lemma~\ref{lem1}
it follows that each cell $2$-representation is simple. So, we have the following claim
(see \cite[Section~3]{MM5}):

\begin{proposition}\label{prop2}
Every  cell $2$-representation of a finitary $2$-category is both simple and transitive.
\end{proposition}

\subsection{Alternative construction of cell $2$-representations for fiat $2$-categories}\label{s4.4}

Let $\cC$ be a fiat $2$-category, $\mathcal{L}$ a left cell in $\cC$ and  
$\mathtt{i}=\mathtt{i}_{\mathcal{L}}$. Consider the principal $2$-representation $\mathbf{P}_{\mathtt{i}}$
and its abelianization $\overline{\mathbf{P}}_{\mathtt{i}}$. For $\mathtt{j}\in\cC$, projective objects
in $\overline{\mathbf{P}}_{\mathtt{i}}(\mathtt{j})$ correspond (up to isomorphism) to $1$-morphisms 
$\mathrm{F}\in\cC(\mathtt{i},\mathtt{j})$ and are denoted $P_{\mathrm{F}}$. The simple top of
$P_{\mathrm{F}}$ is denoted $L_{\mathrm{F}}$.

In the fiat case, \cite[Section~4]{MM1} provides an alternative construction of cell 
$2$-representation which is based on the notion of {\em Duflo involution}.
The following is \cite[Proposition~17]{MM1}:

\begin{proposition}\label{prop3}
The left cell $\mathcal{L}$ contains a unique element 
$\mathrm{G}=\mathrm{G}_{\mathcal{L}}$, called the {\em Duflo involution} in $\mathcal{L}$, 
such that there is a sub-object $K$ of the projective object $P_{\mathbbm{1}_{\mathtt{i}}}$ 
satisfying the following conditions:
\begin{enumerate}[$($a$)$]
\item\label{prop3.1} $\mathrm{F}\, (P_{\mathbbm{1}_{\mathtt{i}}}/K)=0$, for all $\mathrm{F}\in \mathcal{L}$;
\item\label{prop3.2} $\mathrm{F}\, \mathrm{top}(K)\neq 0$, for all $\mathrm{F}\in \mathcal{L}$;
\item\label{prop3.3} $\mathrm{top}(K)\cong L_{\mathrm{G}}$.
\end{enumerate}
\end{proposition}

For $\mathtt{j}\in\cC$, denote by $\mathbf{N}(\mathtt{j})$ the additive closure, in 
$\overline{\mathbf{P}}_{\mathtt{i}}(\mathtt{j})$, of all elements of the form 
$\mathrm{F}\, L_{\mathrm{G}}$, where $\mathrm{F}\in\mathcal{L}\cap\cC(\mathtt{i},\mathtt{j})$.
The following is \cite[Proposition~22]{MM2}:

\begin{proposition}\label{prop4}
{\hspace{2mm}}

\begin{enumerate}[$($i$)$]
\item\label{prop4.1} The assignment $\mathbf{N}$ inherits, by restriction from $\overline{\mathbf{P}}_{\mathtt{i}}$ ,
the natural structure of a $2$-representation of $\cC$.
\item\label{prop4.2} The $2$-representations $\mathbf{N}$ and $\mathbf{C}_{\mathcal{L}}$
are equivalent.
\end{enumerate}
\end{proposition}

The notion of Duflo involution was generalized to some non-fiat $2$-categories in \cite{Zh1}.
In \cite[Example~8]{Xa} (see also \cite[Subsection~9.3]{KMMZ}) one finds an example of a
fiat $2$-category with a left cell $\mathcal{L}$ such that the corresponding Duflo
involution $\mathrm{G}$ satisfies $\mathrm{G}\not\cong\mathrm{G}^{\star}$.

\section{Simple transitive $2$-representations}\label{s5}

\subsection{Weak Jordan-H{\"o}lder theory}\label{s5.1}

In this subsection we overview the weak Jordan-H{\"o}lder theory for additive $2$-representations
of finitary $2$-categories developed in \cite[Section~4]{MM5}. Here 
simple transitive $2$-representations play a crucial role. We start with the following observation
which is just a variation of Lemma~\ref{lem1}, see \cite[Lemma~4]{MM5}.

\begin{lemma}\label{lem11}
Every transitive $2$-representation of a finitary  $2$-category has a unique simple transitive quotient.
\end{lemma}

Let $\cC$ be a finitary $2$-category and $\mathbf{M}\in\cC$-afmod. Let $\mathrm{Ind}(\mathbf{M})$ be the
(finite!) set of isomorphism classes of indecomposable objects in 
\begin{displaymath}
\coprod_{\mathtt{i}\in\ccC} \mathbf{M}(\mathtt{i}).
\end{displaymath}
The {\em action pre-order} $\to_{\ccC}$ on $\mathrm{Ind}(\mathbf{M})$ is defined as follows: 
$X\to_{\ccC}Y$ provided that $Y$ is isomorphic to a direct summand of $\mathrm{F}\, X$, for some
$1$-morphism $\mathrm{F}$ in $\cC$. Consider a filtration
\begin{equation}\label{eq2}
\varnothing=\mathcal{Q}_0\subsetneq \mathcal{Q}_1\subsetneq \dots \subsetneq \mathcal{Q}_m=\mathrm{Ind}(\mathbf{M})
\end{equation}
such that, for each $i=1,2,\dots,m$, 
\begin{itemize}
\item the set $\mathcal{Q}_i\setminus \mathcal{Q}_{i-1}$ is an equivalence class with respect to $\to_{\ccC}$,
\item the set $\mathcal{Q}_i$ has the property that, for all $X\in \mathcal{Q}_i$ and $Y\in \mathrm{Ind}(\mathbf{M})$
such that $X\to_{\ccC}Y$, we have $Y\in \mathcal{Q}_i$.
\end{itemize}
For $i=1,2,\dots,m$ and $\mathtt{j}\in\cC$, let $\mathbf{M}_i(\mathtt{j})$ denote the additive closure in
$\mathbf{M}(\mathtt{j})$ of all objects from $\mathcal{Q}_i\cap \mathbf{M}(\mathtt{j})$. Furthermore, we
denote by $\mathbf{I}_i(\mathtt{j})$ the ideal of $\mathbf{M}_i(\mathtt{j})$ generated by all
objects from $\mathcal{Q}_{i-1}\cap \mathbf{M}(\mathtt{j})$.

The assignment $\mathbf{M}_i$ inherits, by restriction from $\mathbf{M}$, the structure of a $2$-representation
of $\cC$, moreover, $\mathbf{I}_i$ is a $\cC$-stable ideal in $\mathbf{M}_i$. Hence we have a filtration
\begin{displaymath}
0\subset  \mathbf{M}_1\subset  \mathbf{M}_2\subset  \dots\subset \mathbf{M}_m=\mathbf{M}
\end{displaymath}
of $\mathbf{M}$ by $2$-representations. This is a {\em weak Jordan-H{\"o}lder } series of $\mathbf{M}$.
By construction, the $2$-representation $\mathbf{M}_i/\mathbf{I}_i$ is transitive and we can denote by 
$\mathbf{L}_i$ its unique simple transitive quotient, given by Lemma~\ref{lem11}. For a $2$-representation
$\mathbf{N}$, denote by $[\mathbf{N}]$ the equivalence class of $\mathbf{N}$. The multiset
$\{[\mathbf{L}_i]\,:\,i=1,2,\dots,m\}$ is then the multiset of {\em weak Jordan-H{\"o}lder subquotients}
of $\mathbf{M}$. The following is \cite[Theorem~8]{MM5}.

\begin{theorem}\label{thm12}
The multiset  $\{[\mathbf{L}_i]\,:\,i=1,2,\dots,m\}$ does not depend on the choice of the
filtration \eqref{eq2}. 
\end{theorem}

\subsection{How many simple transitive $2$-representations do we have?}\label{s5.2}

Theorem~\ref{thm12} motivates the study of simple transitive $2$-representations for finitary 
$2$-categories. Note that simple $2$-representations are automatically transitive and hence
the name is slightly redundant. However, the name {\em simple transitive} has an advantage that
it addresses both layers of the representation structure:
\begin{itemize}
\item transitivity refers to the discrete layer of objects;
\item simplicity refers to the $\Bbbk$-linear layer of morphisms.
\end{itemize}
For an arbitrary finitary  $2$-category $\cC$, we thus have a general problem:

\begin{problem}\label{prob14}
Classify all simple transitive $2$-representations of $\cC$, up to equivalence. 
\end{problem}

Later on we will survey the cases in which this problem is solved, however,  for general $\cC$,
even for general fiat $\cC$, it is wide open.

Of course, one could try to draw parallels with finite dimensional algebras.
One very easy fact from the classical representation theory is that every finite dimensional 
$\Bbbk$-algebra has only a finite number of isomorphism classes of simple modules. 
The corresponding statement in $2$-representation theory is still open. At the
moment this seems to be one of the major challenges in this theory.

\begin{question}\label{quest15}
Is it true that, for any finitary (fiat) $2$-category $\cC$, the number of equivalence classes of 
simple transitive $2$-representations of $\cC$ is finite?
\end{question}

In all cases in which the answer to Problem~\ref{prob14} is known (see below), the number of 
of equivalence classes of  simple transitive $2$-representations is indeed finite.

\subsection{Finitary $2$-categories of type I}\label{s5.3}

A finitary $2$-category $\cC$ is said to be of {\em type} I provided that every 
simple transitive $2$-representation of $\cC$ is equivalent to a cell $2$-representations.
Thus, for a finitary $2$-category of type I, Problem~\ref{prob14} has a fairly easy answer
(modulo comparison of cell $2$-representations with each other). Moreover, for a finitary 
$2$-category of type I, Question~\ref{quest15} has positive answer. The first example
of finitary $2$-categories of type I is provided by \cite[Theorem~18]{MM5} and \cite[Theorem~33]{MM6}.

\begin{theorem}\label{thm16}
Every weakly fiat $2$-category with strongly regular two-sided cells is of type I. 
\end{theorem}

Here we also note that many results in \cite{MM1}--\cite{MM5} assume that a certain {\em numerical condition}
is satisfied. This assumption was rendered superfluous by \cite[Proposition~1]{MM6}.

As a special case of Theorem~\ref{thm16}, we have that, the $2$-category $\cC_A$,
for a self-injective finite dimensional $\Bbbk$-algebra $A$, see Subsection~\ref{s3.2},
is of type I. Another special case of Theorem~\ref{thm16} is that  the $2$-category $\cS_{W}$
of Soergel bimodules in Weyl type $A$ (that is when $W$ is isomorphic to the symmetric group),
see Subsection~\ref{s3.5}, is of type I. Furthermore, all finitary quotients of 
$2$-Kac-moody algebras from Subsection~\ref{s3.7} are of type I, see \cite[Subsection~7.2]{MM5}
for details. 

Simple transitive $2$-representations of the $2$-category $\cC_A$ were also studied
for some $A$ which are not self-injective (in this case $\cC_A$ is not weakly fiat). 
The first result in this direction was the following statement, which is the main result of \cite{MZ1}.

\begin{theorem}\label{thm18}
For $A=A_2$ or $A_3$ as in Subsection~\ref{s3.3}, the corresponding $2$-category $\cC_A$ is of type I.
\end{theorem}

Despite of the fact that $A$ as in Theorem~\ref{thm18} is not self-injective, it has a non-zero
projective-injective module. The latter plays a crucial role in the arguments. Recently, based on
some progress made in \cite{KM2} and \cite{KMMZ}, Theorem~\ref{thm18} was generalized in \cite{MZ2}
as follows:

\begin{theorem}\label{thm19}
Let $A$  be a basic connected $\Bbbk$-algebra which has a non-zero projective-injective module
and which is directed in the sense that the Gabriel quiver of the algebra $A$ has neither loops nor
oriented cycles. Then the corresponding $2$-category $\cC_A$ is of type I.
\end{theorem}

When $A$ does not have any non-zero projective-injective module, the approach of \cite{MM5,MZ1,MZ2} 
fails. Only one special case (the smallest one) was recently completed in \cite[Theorem~6]{MMZ} (partially based on 
\cite[Subsection~7.1]{MM5}).

\begin{theorem}\label{thm20}
For $A=\Bbbk[x,y]/(x^2,y^2,xy)$, the corresponding $2$-category $\cC_A$ is of type I.
\end{theorem}

After all the cases listed above, the following question is rather natural:

\begin{question}\label{quest17}
Is it true that the $2$-category $\cC_A$ is of type I, for any $A$?
\end{question}

Apart from the cases listed above, there is a number of other type I examples.
The following result in \cite[Theorem~6.1]{Zi}.

\begin{theorem}\label{thm21}
The $2$-category of Soergel bimodules in Weyl type $B_2$ is of type I.
\end{theorem}

One interesting difference of the latter case compared to all other cases listed above
is the fact that, for the $2$-category of Soergel bimodules in Weyl type $B_2$, there are
non-equivalent cell $2$-representations which correspond to left cells inside the same two-sided cell.

Let $(W,S)$ be a finite Coxeter system. The corresponding $2$-category $\cS_W$ of Soergel bimodules,
see Subsection~\ref{s3.5}, has a unique minimum two-sided cell consisting of the identity $1$-morphism.
If we take this minimum two-sided cell away, in what remains there is again a unique minimum two-sided
cell $\mathcal{J}$. This two-sided cell contains, in particular, all Soergel bimodules of the form
$\mathtt{C}\otimes_{\mathtt{C}^s}\mathtt{C}$, where $s\in S$. There is a unique $2$-ideal $\cI$ in 
$\cS_W$ which is maximal, with respect to inclusions, in the set of all $2$-ideals in 
$\cS_W$ that do not contain identity $2$-morphisms for $1$-morphisms in $\mathcal{J}$. The quotient
$\underline{\cS_W}:=\cS_W/\cJ$ is called the {\em small quotient} of $\cS_W$, see \cite[Subsection~3.2]{KMMZ}.
The $2$-category $\underline{\cS_W}$ inherits from $\cS_W$ the structure of a fiat $2$-category.
The following result can be found in \cite[Sections~6, 7 and 8]{KMMZ}.

\begin{theorem}\label{thm22}
{\hspace{2mm}}

\begin{enumerate}[$($i$)$]
\item\label{thm22.1} If $|S|>2$, then  $\underline{\cS_W}$ is of type I.
\item\label{thm22.2} If $(W,S)$ is of Coxeter type $I_2(n)$, with $n>4$, then 
$\underline{\cS_W}$ is of type I if and only if $n$ is odd.
\end{enumerate} 
\end{theorem}

\subsection{Finitary $2$-categories that are not of type I}\label{s5.4}

First, rather degenerate examples of finitary $2$-categories which are not of type I were constructed
already in \cite[Subsection~3.2]{MM5}. They are inspired by transitive group actions. Each finite group 
structure can be extended to a fiat $2$-category in a fairly obvious way (by adding formal direct sums
of elements and linearizing spaces of identity $2$-morphisms). The resulting $2$-category has just one left cell 
and the corresponding cell $2$-representation is, morally, the left regular representation of the group.
However, simple transitive $2$-representations correspond to transitive actions of the original group
on sets. The latter are given by action on (left) cosets modulo subgroups. In particular, we get a lot of
simple transitive $2$-representations which are not cell $2$-representations. The example, however, feels
rather artificial.

The first more ``natural'' example was constructed in \cite{MaMa}. Let $(W,S)$ be a Coxeter system
of type $I_2(n)$, with $n>3$. Consider the small quotient $\underline{\cS_W}$ of the
$2$-category of Soergel bimodules. This is a fiat $2$-category with two two-sided cells. The minimum
one consists just of the identity $1$-morphisms. The maximum one is not strongly regular. Let $S=\{s,t\}$.
Denote by $\cQ_n$ the $2$-full sub-$2$-category of $\underline{\cS_W}$ given by all $1$-morphisms in the additive
closure of the identity $1$-morphisms and of all $1$-morphisms $\mathrm{F}$ which lie in the same right
cell and in the same left cell as the $1$-morphisms given by $\mathtt{C}\otimes_{\mathtt{C}^s}\mathtt{C}$.
The main result of \cite{MaMa} is the following:

\begin{theorem}\label{thm23}
{\hspace{2mm}}

\begin{enumerate}[$($i$)$]
\item\label{thm23.1} The $2$-category  $\cQ_5$ is of type I.
\item\label{thm23.2} The $2$-category  $\cQ_4$ is not of type I. In fact, $\cQ_4$ has a unique 
(up to equivalence) simple transitive $2$-representation which is not equivalent to any of two
cell $2$-representations. 
\end{enumerate}
\end{theorem}

The major part of \cite{MaMa} is devoted to an explicit construction of this additional 
simple transitive $2$-representation of $\cQ_4$. The construction is very technical and is based
on the following idea: One of the two cell $2$-representations of $\cQ_4$ has an invertible
automorphism which swaps the isomorphism classes of the two non-isomorphic
indecomposable objects in the underlying
category of this $2$-representation. The additional
simple transitive $2$-representation is constructed using the orbit category with respect to this
non-trivial automorphism. The main issue is that this automorphism is not strict (as homomorphisms
between $2$-rep\-re\-sen\-ta\-tions are not strict) and so it requires a lot of technical effort to
go around this complication. Based on this construction, the following statement was proved
in \cite[Sections~7]{KMMZ}.

\begin{theorem}\label{thm24}
Let $(W,S)$ be of Coxeter type $I_2(n)$, with $n>4$ even. Then $\underline{\cS_W}$ is not of type I.
Moreover, the following holds:
\begin{enumerate}[$($i$)$]
\item\label{thm24.1} Two of the three cell $2$-representations of $\underline{\cS_W}$ have an invertible
automorphism which is not isomorphic to the identity. The orbit construction as in \cite{MaMa} 
with respect to this automorphism produces a new simple transitive $2$-representation of $\underline{\cS_W}$.
\item\label{thm24.2} If $n\neq 12,18,30$, then every 
simple transitive $2$-representation of $\underline{\cS_W}$ is equivalent to either a cell $2$-representation
or one of the $2$-representations constructed in \eqref{thm24.1}.
\end{enumerate}
\end{theorem}

\subsection{Schur's lemma}\label{s5.5}
As we saw in the previous subsection, endomorphisms of cell $2$-representations play an important role
in this study. This naturally raises the following problem:

\begin{problem}\label{quest25}
Describe the (properties of the) bicategory of endomorphisms  of a simple transitive 
$2$-representation of a finitary $2$-category.
\end{problem}

The only known result in this direction is the following statement which is \cite[Theorem~16]{MM3}.

\begin{theorem}\label{thm26}
Let $\cC$ be a fiat $2$-category, $\mathcal{J}$ a strongly regular two-sided cell in $\cC$
and $\mathcal{L}$ a left cell in $\mathcal{J}$. Then any endomorphism of the cell
$2$-representation $\mathbf{C}_{\mathcal{L}}$ is isomorphic to the direct sum of a number
of copies of the identity endomorphism $\mathrm{ID}_{\mathbf{C}_{\mathcal{L}}}$. 
Moreover, the endomorphism space (given by all modifications)
of $\mathrm{ID}_{\mathbf{C}_{\mathcal{L}}}$ consist just of scalar multiplies of the identity
modification.
\end{theorem}

\subsection{Apex}\label{s5.6}

The following is \cite[Lemma~1]{CM}.

\begin{lemma}\label{lem31}
Let $\cC$ be a finitary $2$-category and $\mathbf{M}$ a transitive $2$-representation of $\cC$.
There is a unique two sided cell $\mathcal{J}=\mathcal{J}_{\mathbf{M}}$ which is maximal, with respect 
to the two-sided order,  in the set of all two-sided cells that contain $1$-morphisms which are not
annihilated by $\mathbf{M}$. The two-sided cell $\mathcal{J}$ is idempotent.
\end{lemma}

The two-sided cell $\mathcal{J}$ is called the {\em apex} of $\mathbf{M}$. The general problem of 
classification of all simple transitive $2$-representations of $\cC$ thus splits naturally into 
subproblems to classify simple transitive $2$-representations of $\cC$ with a fixed apex
$\mathcal{J}$ (which should be an idempotent two-sided cell). Quite often, this simplifies the problem,
due to the following result which is proved analogously to \cite[Theorem~18]{MM5}.

\begin{theorem}\label{thm31}
Let $\cC$ be a weakly fiat $2$-category, $\mathcal{J}$ a strongly regular two-sided cell in $\cC$ and
$\mathbf{M}$ a simple transitive $2$-representation of $\cC$ with apex $\mathcal{J}$.
Then $\mathbf{M}$ is equivalent to a cell $2$-representation.
\end{theorem}

\subsection{Connection to integral matrices}\label{s5.7}

Let $\cC$ be a finitary $2$-category and $\mathbf{M}\in\cC$-afmod. Let $\mathrm{Ind}(\mathbf{M})$ be
as in Subsection~\ref{s5.1}. Then, to any $1$-morphism $\mathrm{F}$ in $\cC$, one can associate a matrix
$[\mathrm{F}]_{\mathbf{M}}$ whose rows and columns are indexed by elements in $\mathrm{Ind}(\mathbf{M})$
and the intersection of the $X$-row and $Y$-column gives the multiplicity of $X$ as a summand of
$\mathrm{F}\, Y$. The following observation is \cite[Lemma~11(ii)]{MM5}.

\begin{lemma}\label{lem32}
Assume that $\mathbf{M}$ is transitive and that  $\mathrm{F}$ contains, as summands, representatives
from all isomorphism classes of indecomposable $1$-morphisms in $\cC$. Then 
all coefficients in $[\mathrm{F}]_{\mathbf{M}}$ are positive.
\end{lemma}

This observation allows one to use the classical Perron-Frobenius Theorem  (see \cite{Pe,Fro1,Fro2}) in 
the study of  simple transitive $2$-representations. This is an important ingredient in the arguments
in \cite{MM5,MZ1,Zi,MMZ}. 

The above observation also provides some evidence for the general positive answer to 
Question~\ref{quest15}. 
Indeed, the Grothendieck decategorification of a finitary $2$-category $\cC$
gives a finite dimensional $\Bbbk$-algebra, call it $A$. For each $1$-morphism $\mathrm{F}$ in $\cC$, 
we thus have the minimal polynomial $g_{\mathrm{F}}(\lambda)$ for the class $[\mathrm{F}]_{\oplus}$ in $A$. 
Now, if $\mathbf{M}\in\cC$-afmod, then $[\mathbf{M}]_{\oplus}$ gives rise to an $A$-module and hence
\begin{displaymath}
g_{\mathrm{F}}([\mathrm{F}]_{\mathbf{M}})=0.
\end{displaymath}
Therefore, if $\mathbf{M}$ is transitive and $\mathrm{F}$ contains, as summands, representatives
from all isomorphism classes of indecomposable $1$-morphisms in $\cC$, then, because of 
Lemma~\ref{lem32} and \cite[Theorem~3.2]{Zi2}, there are only finitely many possibilities
for the matrix $[\mathrm{F}]_{\mathbf{M}}$. Consequently, there are only finitely many possibilities 
for matrices $[\mathrm{G}]_{\mathbf{M}}$, where $\mathrm{G}$ is an indecomposable $1$-morphism in $\cC$.

In many papers, for instance, in  \cite{MM5,MM6,MZ1,MZ2,Zi}, the classification problem was approached
in two steps. The first step addressed classification of all possibilities for matrices 
$[\mathrm{F}]_{\mathbf{M}}$. The second step studied actual $2$-representations for each solution
provided by  the first step. In case of the $2$-category $\cC_A$, the first step studied matrices
$M$ with positive integer coefficients satisfying $M^2=\dim(A) M$. This is an interesting combinatorial
problem which was investigated in detail in \cite{Zi2}.

\subsection{Approach using (co)algebra objects}\label{s5.8}

The story with the cases $n=12$, $18$, $30$ in Theorem~\ref{thm24}\eqref{thm24.2} was quite interesting.
The detailed study of integral matrices, as outlined in Subsection~\ref{s5.7}, suggested in these cases 
possibility of existence of simple transitive $2$-representations of $\cS_W$ which are neither 
cell $2$-representations nor the ones constructed in Theorem~\ref{thm24}\eqref{thm24.1}.
These $2$-representations were constructed in \cite{MT}, based on \cite{El}, using diagrammatic calculus.
Under the additional assumption of gradeability, it was shown that, together with 
cell $2$-representations and the $2$-representations constructed in Theorem~\ref{thm24}\eqref{thm24.1},
these exhaust all simple transitive $2$-representations of $\cS_W$.

Unfortunately, the diagrammatic calculus is not really compatible with our definitions of $2$-categories.
This raised a natural problem to reformulate the results of \cite{MT} in some language compatible
with our definitions. This was achieved in \cite{MMMT} using ideas of \cite{EGNO} related to the study of
algebra and coalgebra objects in $2$-categories. More precisely, the following is \cite[Theorem~9]{MMMT}.

\begin{theorem}\label{thm33}
Let $\cC$ be a fiat $2$-category and $\mathbf{M}$ a transitive $2$-representations of $\cC$. Then there
is a coalgebra object $A$ in the injective abelianization $\underline{\cC}$ of $\cC$ such that 
$\mathbf{M}$ is equivalent  to the $2$-representation of $\cC$ given by the action of $\cC$ on the 
category of injective right $A$-comodule objects in $\underline{\cC}$.
\end{theorem}

For cell $2$-representations, the corresponding coalgebra objects turn out to be related to 
Duflo involutions, see \cite[Subsection~6.3]{MMMT}.
Theorem~\ref{thm33} motivates the following general problem:

\begin{problem}\label{prob34}
Classify all coalgebra objects in $\underline{\cC}$, up to isomorphism. 
\end{problem}

In this general formulation, Problem~\ref{prob34} is certainly more difficult than Problem~\ref{prob14}.
However, the very useful side of Theorem~\ref{thm33} is that one can construct simple transitive 
$2$-representations by guessing the corresponding coalgebra objects (as it was done for the ``exotic''
simple transitive $2$-representations of $\cS_W$ in types $I_2(12)$, $I_2(18)$ and $I_2(30)$ in
\cite{MMMT}).

\section{Other classes of $2$-representations and related questions}\label{s6}

\subsection{Isotypic $2$-representations}\label{s6.1}

Let $\cC$ be a finitary $2$-category and $\mathbf{M}$ a $2$-rep\-re\-sen\-ta\-tion of $\cC$.
We will say that $\mathbf{M}$ is {\em isotypic} provided that all weak Jordan-H{\"o}lder
subquotients of $\mathbf{M}$ are equivalent, see \cite[Subsection~4.3.4]{Ro2} and \cite[Subsection~3.6]{MM6}.

For any $2$-representation $\mathbf{M}$ of $\cC$ and any finitary $\Bbbk$-linear 
$2$-category $\mathcal{A}$ one defines the {\em inflation} $\mathbf{M}^{\boxtimes\mathcal{A}}$
of $\mathbf{M}$ by $\mathcal{A}$ as the $2$-representation of $\cC$ which sends each
$\mathtt{i}\in \cC$ to the tensor product $\mathbf{M}(\mathtt{i})\boxtimes \mathcal{A}$
and defines the action of $\cC$ on the objects and morphisms in these tensor products
by acting on the first component, see \cite[Subsection~3.6]{MM6} for details.
The following result is \cite[Theorem~4]{MM6}.

\begin{theorem}\label{thm37}
Let $\cC$ be a weakly fiat $2$-category with a unique maximal two-sided cell 
$\mathcal{J}$. Let $\mathcal{L}$ be a left cell in $\mathcal{J}$. Assume that 
$\mathcal{J}$ is strongly regular and that any non-zero $2$-ideal of $\cC$ contains the identity 
$2$-morphism, for some $1$-morphism in $\mathcal{J}$.
Then any isotypic faithful $2$-representation of $\cC$ is equivalent to an inflation of 
the cell $2$-representation $\mathbf{C}_{\mathcal{L}}$. 
\end{theorem}

For finitary quotients of $2$-Kac-Moody algebras, the statement of Theorem~\ref{thm37}
can be deduced from \cite[Subsection~4.3.4]{Ro2}. Compared to the general case, the case of
$2$-Kac-Moody algebras is substantially simplified by existence of idempotent $1$-morphisms
in each two-sided cell. A challenging problem related to isotypic $2$-representation is the following:

\begin{problem}\label{thm37}
Classify faithful isotypic $2$-representation for an arbitrary 
weakly fiat $2$-category $\cC$ with unique maximal two-sided cell $\mathcal{J}$
under the assumption that that $\mathcal{J}$ is strongly regular.
\end{problem}

In the easiest case, this problem will appear in the next subsection.

\subsection{All $2$-representations}\label{s6.2}

The question of classification of all $2$-representations, for a given finitary $2$-category
$\cC$, is open in all non-trivial case. The only {\em trivial} case is the case when the only indecomposable
$1$-morphisms in $\cC$ are the identities, up to isomorphism. It is certainly enough to consider
the case when $\cC$ has one object, say $\mathtt{i}$. Up to biequivalence, we may also assume
that $\mathbb{1}_{\mathtt{i}}$ is the only indecomposable $1$-morphism in $\cC$
(on the nose and not just up to isomorphism). Then, directly from the definition,
we have that a  $2$-representation of such $\cC$ is given by a pair $(Q,\varphi)$,
where $Q$ is a finite dimensional $\Bbbk$-algebra and $\varphi$ is an algebra homomorphism from
$\mathrm{End}(\mathbb{1}_{\mathtt{i}})$ to the center of $Q$. 
Then $\cC$ acts on a small category equivalent to
$Q$-proj in the obvious way. Note that all $2$-representations
of $\cC$ are isotypic. Furthermore, $\cC$ satisfies the assumptions of Theorem~\ref{thm37}
if and only if $\mathrm{End}(\mathbb{1}_{\mathtt{i}})\cong\Bbbk$.

The first non-trivial case to consider would be the following:

\begin{problem}\label{prop45}
Classify, up to equivalence, all finitary additive $2$-representations of $\cC_D$,
where $D=\mathbb{C}[x]/(x^2)$.
\end{problem}

\subsection{Discrete extensions between $2$-representations}\label{s6.3}

A major challenge in $2$-representation theory is the following:

\begin{problem}\label{prob46}
Develop a sensible homological theory for the study of $2$-rep\-re\-sen\-ta\-tions.
\end{problem}

A fairly naive attempt to define some analogue of $\mathrm{Ext}^1$ for $2$-representations
was made in \cite{CM}. 

Let $\cC$ be a finitary $2$-category and $\mathbf{M}\in\cC$-afmod. For each $\mathtt{i}\in\cC$,
choose a full, additive, idempotents split and isomorphism closed subcategory $\mathbf{K}(\mathtt{i})$
of $\mathbf{M}(\mathtt{i})$ such that $\mathbf{K}$ becomes a sub-$2$-representation of $\cC$
by restriction. Let $\mathbf{I}$ be the ideal of $\mathbf{M}$ generated by $\mathbf{K}$
and $\mathbf{N}:=\mathbf{M}/\mathbf{I}$. Then the sequence
\begin{equation}\label{eq3}
0\to \mathbf{K}\overset{\Phi}{\longrightarrow}\mathbf{M}\overset{\Psi}{\longrightarrow}\mathbf{N}\to 0,
\end{equation}
where $\Phi$ is the natural inclusion and $\Psi$ is the natural projection, will be called a
{\em short exact sequence} of $2$-representations of $\cC$. The {\em discrete extension} 
$\Theta$ {\em realized} by \eqref{eq3} is the subset of $\mathcal{S}[\cC]$ that consists of all
classes $[\mathrm{F}]$ for which there exist an indecomposable object $X$ in some
$\mathbf{M}(\mathtt{i})\setminus \mathbf{K}(\mathtt{i})$ and an indecomposable object $Y$ in some
$\mathbf{K}(\mathtt{j})$ such that $Y$ is isomorphic to a summand of $\mathrm{F}\, X$.

For $\mathbf{K}',\mathbf{N}'\in\cC$-afmod, the set $\mathrm{Dext}(\mathbf{N}',\mathbf{K}')$ of
{\em discrete extensions} from $\mathbf{N}'$ to $\mathbf{K}'$
consists of all possible $\Theta$ which are realized by some short
exact sequence \eqref{eq3} with $\mathbf{K}$ equivalent to $\mathbf{K}'$ and 
$\mathbf{N}$ equivalent to $\mathbf{N}'$.

In many case, a very useful piece of information is to know whether 
$\mathrm{Dext}(\mathbf{N},\mathbf{K})$ is empty (i.e. the first discrete extension {\em vanishes}) 
or not. In a number of cases, one could also either explicitly describe all elements in 
$\mathrm{Dext}(\mathbf{N},\mathbf{K})$ or at least give a reasonable estimate of how they look like.
Vanishing of discrete extensions between certain simple transitive $2$-representations appears
in a disguised form and is an essential part of the arguments in \cite{MM5,MM6}.
The following is \cite[Theorem~25]{CM}.

\begin{theorem}\label{thm47}
Let $\cC$ be a weakly fiat $2$-category, 
$\mathbf{K}$ a transitive $2$-representation of $\cC$ with apex $\mathcal{J}_{\mathbf{K}}$,
and $\mathbf{N}$ a transitive $2$-representation of $\cC$ with apex $\mathcal{J}_{\mathbf{N}}$.
Assume that, for any left cell $\mathcal{L}$ in $\mathcal{J}_{\mathbf{N}}$,
there exists a left cell $\mathcal{L}'$ in $\mathcal{J}_{\mathbf{K}}$ such that 
$\mathcal{L}\geq_L\mathcal{L}'$. Then $\mathrm{Dext}(\mathbf{N},\mathbf{K})=\varnothing$.
\end{theorem}

As a consequence of Theorem~\ref{thm47}, all discrete self-extensions for transitive 
$2$-rep\-re\-sen\-ta\-tions of weakly fiat $2$-categories vanish.

The results in \cite[Subsection~7.2]{CM} suggest that the answer to Question~\ref{quest48}
might be interesting and is not obvious.

\begin{question}\label{quest48}
What is $\mathrm{Dext}(\mathbf{N},\mathbf{K})$, for any pair 
$(\mathbf{N},\mathbf{K})$ of simple transitive $2$-representations of the
$2$-category $\cS_W$ of Soergel bimodules in Weyl type $A$?
\end{question}

\subsection{Applications}\label{s6.4}

The first, rather spectacular, application of classification of certain classes of $2$-representations
appears in \cite{CR}. More precisely, \cite[Proposition~5.26]{CR} classifies, up to equivalence, 
all $2$-representations of the $2$-Kac-Moody version of $\mathfrak{sl}_2$ which satisfy a number 
of natural assumptions. This is an essential ingredient in the proof of derived equivalence for certain
blocks of the symmetric group, see  \cite[Theorem~7.6]{CR}. Similar ideology was used, in particular, 
to describe blocks of Lie superalgebras, see, for example, \cite{BS,BLW} and references therein.

In \cite{KM1}, classification of simple transitive $2$-representation for the $2$-category of
Soergel bimodules in type $A$ (cf. Theorem~\ref{thm16}) was used to classify indecomposable projective functors on the
principal block of BGG category $\mathcal{O}$ for $\mathfrak{sl}_n$.

\vspace{2mm}

\noindent
V.~M.: Department of Mathematics, Uppsala University, Box. 480,
SE-75106, Uppsala, SWEDEN, email: {\tt mazor\symbol{64}math.uu.se}

\end{document}